\documentclass{article}
\usepackage{amsmath}
\usepackage{amssymb}
\usepackage{bbm}
\usepackage{amsthm}
\usepackage{enumitem}
\usepackage{authblk}
\usepackage{mathtools}
\usepackage[margin=1.2in]{geometry}
\usepackage[backend=biber, style=trad-plain, sorting=nty]{biblatex}
\begin{document}
\title {Short proof of the sharpness of the phase transition for the random-cluster model with $q=2$  }
\author{Yacine Aoun}
\affil{Université de Genève}
\maketitle
         Abstract: The purpose of this modest note is to provide a short proof of the sharpness of the phase transition for the Random-cluster model with $q=2$ by extending the approach developed by Duminil-Copin and Tassion \cite{HT} for $q=1$. This in particular implies the exponential decay of the two point-correlation function in the subcritical Ising model.         
\section{Introduction}        
Let us start by defining the nearest neighbor random-cluster measure on $\mathbb{Z}^{d}$. For a finite subgraph $\Lambda=(V,E)$ of $\mathbb{Z}^{d}$, a percolation configuration $\omega=(\omega)_{e\in E}$ is an element of $\lbrace 0,1\rbrace^{E}$. A configuration $\omega$ can be seen as a subgraph of $\Lambda$ with vertex-set $V$ and edge-set given by $\lbrace \lbrace x,y\rbrace\in E : \omega_{x,y}=1\rbrace$. If $\omega_{x,y}=1$, we say that $\lbrace x,y\rbrace$ is open. Let $k(\omega)$ be the number of connected components in $\omega$ and $o(\omega)$ (respectively $f(\omega)$) the number of open (respectively closed) edges in $\omega$.

Fix $p\in [0,1], q>0$. Let $\mu_{\Lambda,p,q}$ be a measure defined for any $\omega\in\lbrace 0,1\rbrace^{E}$ by 
\begin{center}
$\mu_{\Lambda,p,q}(\omega)=\dfrac{q^{k(\omega)}}{Z}p^{o(\omega)}(1-p)^{f(\omega)}$,
\end{center}
where $Z$ is a normalizing constant introduced in such a way that $\mu_{\Lambda,p,q}$ is a probability measure. The measure $\mu_{\Lambda,p,q}$ is called the random-cluster measure on $\Lambda$ with free boundary conditions. For $q\geq 1$, the measures can be extented to $\mathbb{Z}^{d}$ by taking the weak limit of measures defined in finite volume. 

We say that $x$ and $y$ are connected in $S\subseteq \mathbb{Z}^{d}$ if there exists a finite sequence of vertices $(v_i)_{i=0}^{n}$ in $S$ such that $v_0 =x$, $v_n =y$ and $ \lbrace v_{i} ,v_{i+1}\rbrace$ is \textit{open} for every $0\leq i<n$. We denote this event by $x \overset{S}{\leftrightarrow} y$. For $A\subset\mathbb{Z}^{d}$, we write $x\overset{S}{\leftrightarrow} A$ for the event that $x$ is connected in $S$ to a vertex in $A$. If $S=\mathbb{Z}^{d}$, we drop it from the notation. We write $0\leftrightarrow\infty$ if for every $n\in\mathbb{N}$, there exists $x\in\mathbb{Z}^{d}$ such that $0\leftrightarrow x$ and $\vert x\vert\geq n$, where $\vert \cdot\vert$ denotes a norm on $\mathbb{Z}^{d}$.

For $q\geq 1$, the model undergoes a phase transition: there exists $p_{c}\in[0,1]$ satisfying
\begin{center}
$\mu_{\mathbb{Z}^{d},p,q}(0\leftrightarrow \infty)=
\begin{cases} 
=0 & \text{if } p<p_{c}, \\
>0 &  \text{if } p>p_{c}. \\
\end{cases}
$
\end{center}
A very nice idea introduced in \cite{HT} is to define a new critical parameter $\tilde{p}_{c}$ for which it is easier to prove sharpness, which will in turn imply that $p_{c}=\tilde{p}_{c}$ (see Theorem 1 below). For a finite subset S of $\mathbb{Z}^{d}$, let $\Delta S$ be the set of edges with exactly one endpoint in $S$ and define 
\begin{equation}
\varphi_ p(S) := p \sum\limits_{\lbrace x,y\rbrace\in  \Delta S}\mu_{S,p,q}(0\leftrightarrow x).
\end{equation}
Then define the following critical parameter
\begin{center}
$\tilde{p}_{c}:=\sup\lbrace p\in [0,1] : \varphi_{p}(S)<1$ for some finite $S \subset \mathbb{Z}^{d}$ containing $0\rbrace$.
\end{center}
The main theorem of this note is the following one.
\newtheorem{theorem}{Theorem}
\begin{theorem}
\begin{enumerate}
\item For $p>\tilde{p}_{c}$, $\mu_{\mathbb{Z}^{d},p,2}(0\leftrightarrow\infty)\geq\dfrac{p-\tilde{p}_{c}}{p} $.
\item For $p<\tilde{p}_{c}$, there exists $c=c(p)>0$ such that for every $x\in\mathbb{Z}^{d}$
\begin{equation}
\mu_{\mathbb{Z}^{d},p,2}(0\leftrightarrow x)\leq\exp(-c\vert x\vert).
\end{equation}
\end{enumerate}
\end{theorem}
\newtheorem{corollary}{Corollary}
\begin{corollary}
We have $p_{c}=\tilde{p}_{c}$. In particular, the phase transition for the random-cluster model with $q=2$ is sharp.
\end{corollary}
\begin{corollary}
In the Ising model, the two-point correlation function decays exponentially fast with distance.
\end{corollary}
\noindent
Corollary 1 follows directly from Theorem 1. Corollary 2 follows from the Edward-Sokal coupling (see \cite{GRIM2}).
\section{Proof of Theorem 1}
We will write $\mu_{\Lambda,p}$ instead of $\mu_{\Lambda,p,2}$. Let us start by proving the second item of Theorem 1. Firstly, we will need the following lemma.
\newtheorem{lemma}{Lemma}
\begin{lemma}[Modified Simon's inequality]
Let S be a finite set of $\hphantom{,}\mathbb{Z}^{d}$ containing $0$. For every $z\notin S$,
\begin{equation}
\mu_{\mathbb{Z}^{d},p}(0\leftrightarrow z)\leq p\sum\limits_{\lbrace x,y\rbrace\in\Delta S}\mu_{S,p}(0\leftrightarrow x) \mu_{\mathbb{Z}^{d},p}(y\leftrightarrow z).
\end{equation}
\end{lemma}
A similar inequality was proved for the Ising model in \cite{HT}. Lemma 1 follows from the latter by the Edward-Sokal coupling by remarking that $\tanh(-\frac{1}{2}\log(1-p))\leq p $.

Fix $p<\tilde{p}_{c}$ and $S$ a finite set containing 0 such that $\varphi_{p}(S)<1$. Let $\Lambda_{n}$ be the box of size $n$ around $0$ for the norm $\vert\cdot\vert$. Fix $\Lambda_{L}$ such that $S\subset\Lambda_{L}$. Then, using Lemma 1, we can write
\begin{equation}
\mu_{\mathbb{Z}^{d},p}(0\leftrightarrow z)\leq p\sum\limits_{\lbrace x,y\rbrace\in\Delta S}\mu_{S,p}(0\leftrightarrow x) \mu_{\mathbb{Z}^{d},p}(y\leftrightarrow z)\leq\varphi_{p}(S)\max\limits_{y\in\Lambda_{L}}\mu_{\mathbb{Z}^{d},p}(y\leftrightarrow z).
\end{equation}
Note that $\vert y-z\vert\geq \vert z\vert-L$. If $\vert y-z\vert\leq L$, we bound $\mu_{\mathbb{Z}^{d},p}(y\leftrightarrow z)$ by $1$, otherwise we apply $(4)$ to $y$ and $z$ instead of $0$ and $z$. Iterating $\lfloor\vert z\vert/L\rfloor$ this strategy yields
\begin{center}
$\mu_{\mathbb{Z}^{d},p}(0\leftrightarrow z)\leq\varphi_{p}(S)^{\lfloor \vert z\vert/L\rfloor}$,
\end{center}
which proves the second item of Theorem 1.

We now turn to the proof of the first item of Theorem 1.  Let $p>\tilde{p}_{c}$ and $\partial\Lambda_{n}$ be the boundary of $\Lambda_{n}$. We will prove the following differential inequality. 
\begin{lemma}
Fix $p>\tilde{p}_{c}$. Then
\begin{equation}
\frac{d}{dp}\mu_{\Lambda_{n},p}(0\leftrightarrow \partial\Lambda_{n})\geq\frac{1}{p}(1-\mu_{\Lambda_{n},p}(0\leftrightarrow\partial\Lambda_{n})).
\end{equation}
\end{lemma}
Integrating this inequality between $\tilde{p}_{c}$ and $p$ and taking $n$ to infinity yields the first item of Theorem 1. We will therefore focus on proving Lemma 2. Let $E(\Lambda_{n})$ be the set of edges whose endpoints are in $\Lambda_{n}$.  We will need the following result. 
\begin{lemma}
Let $A$ be an increasing event depending on edges of $\Lambda_{n}$ only. Then 
\begin{equation}
\dfrac{d}{dp}\mu_{\Lambda_{n},p,q}(A)=\sum_{e\in E(\Lambda_{n})}\mu_{\Lambda_{n},p,q}(A \vert\hphantom{,}\omega_{e}=1)-\mu_{\Lambda_{n},p,q}(A \vert\hphantom{,}\omega_{e}=0).
\end{equation}
\end{lemma}
\noindent 
The proof is a straightforward computation. Recall that an edge $e$ is pivotal for a configuration $\omega$ and an event $A$ if $\omega_{(e)}\notin A$ and $\omega^{(e)}\in A$, where $\omega_{(e)}$ (respectively $\omega^{(e)}$) is the same configuration as $\omega$ except maybe for $e$ where we close the edge $e$ in $\omega_{(e)}$ (respectively open the edge $e$ in $\omega^{(e)}$). We can use Lemma 3 and the FKG inequality to see that 
\begin{align*}
\frac{d}{dp}\mu_{\Lambda_{n},p}(A)
&=
\sum\limits_{e\in E(\Lambda_{n})}\mu_{\Lambda_{n},p}(A\vert\omega_{e}=1)-\mu_{\Lambda_{n},p}(A\vert\omega_{e}=0)
\\
&\geq 
\sum\limits_{e\in E(\Lambda_{n})}\mu_{\Lambda_{n},p}(\omega^{(e)}\in A)-\mu_{\Lambda_{n},p}(\omega_{(e)}\in A)
\\
&=
\sum\limits_{e\in E(\Lambda_{n})}\mu_{\Lambda_{n},p}(\text{e pivotal for A}).
\end{align*}
Set $A:=\lbrace 0\leftrightarrow\partial\Lambda_{n}\rbrace$. Define the following random set
\begin{align*}
\gamma:=\lbrace z\in\Lambda_{n} : z\hphantom{:} \text{not connected to}\hphantom{:}\Lambda_{n}^{c}\rbrace .
\end{align*}
By inclusion of events, we get
\begin{align*}
\sum\limits_{e\in E(\Lambda_{n})}\mu_{\Lambda_{n},p}(\text{e pivotal for}\hphantom{:} 0\leftrightarrow\partial\Lambda_{n})
&\geq
\sum\limits_{e\in E(\Lambda_{n})}\mu_{\Lambda_{n},p}(\text{e pivotal for}\hphantom{:} 0\leftrightarrow\partial\Lambda_{n}, 0\nleftrightarrow \partial\Lambda_{n})
\\
&=
\sum\limits_{\substack{S \\ 0\in S}}\sum\limits_{e\in E(\Lambda_{n})}\mu_{\Lambda_{n},p}(\text{e pivotal for}\hphantom{:} 0\leftrightarrow\partial\Lambda_{n}, \gamma=S),
\end{align*}
where we decomposed with respect to all possibilities for $\gamma$ in the last line. Remark that $\gamma =S$ and $e=xy$ is pivotal for $0\leftrightarrow\partial\Lambda_{n}$ if and only if $\gamma=S$, $0\overset{S}{\leftrightarrow}x$ and $y\notin S$. Moreover, the event $\lbrace 0\overset{S}{\leftrightarrow}x\rbrace$ is mesurable with respect to the edges in $S$ and the event $\lbrace\gamma=S\rbrace$ is mesurable with respect to the edges that have at least one endpoint outside of $S$. Finally, all the edges in $\Delta S$ are closed. Thus
\begin{center}
$\mu_{\Lambda_{n},p}(\text{e pivotal for}\hphantom{:} 0\leftrightarrow\partial\Lambda_{n}, \gamma=S)=\mu_{\Lambda_{n},p}(0\overset{S}{\leftrightarrow}x,\gamma=S)=\mu_{S,p}(0\leftrightarrow x)\mu_{\Lambda_{n},p}(\gamma=S)$,
\end{center}  
where the last equality follows from the Markov property. Plugging this into the inequality above gives
\begin{align*}
\sum\limits_{\substack{S \\ 0\in S}}\sum\limits_{e\in E(\Lambda_{n})}\mu_{\Lambda_{n},p}(\text{e pivotal for}\hphantom{:} 0\leftrightarrow\partial\Lambda_{n}, \gamma=S)
&=
\frac{1}{p}\sum\limits_{\substack{S \\ 0\in S}}\sum\limits_{xy\in \Delta S}p\mu_{S,p}(0\leftrightarrow x)\mu_{\Lambda_{n},p}(\gamma=S)
\\
&= \frac{1}{p}\sum\limits_{\substack{S \\ 0\in S}}\varphi_{p}(S)\mu_{\Lambda_{n},p}(\gamma=S)
\\
&\geq 
\frac{1}{p}\mu_{\Lambda_{n},p}(0\nleftrightarrow\partial\Lambda_{n}),
\end{align*}
where we used that $\varphi_{p}(S)\geq 1$ since $p>\tilde{p}_{c}$. Therefore, by combining all the inequalities, we get (5), which finishes the proof of Lemma 2.
\section{Concluding remarks}
1. It is natural to ask whether this approach can be further generalized for bigger values of $q$. The proof of (5) does not use the fact that $q=2$ and is valid for all $q\geq 1$, which implies $\tilde{p}_{c}(q)\geq p_{c}(q)$. However, 
it is easy to see that the susceptibility with free boundary conditions is always infinite at $\tilde{p}_{c}$, i.e. 
\begin{center}
$\sum\limits_{x\in\mathbb{Z}^{d}}\mu_{\mathbb{Z}^{d},\tilde{p}_{c},q}(0\leftrightarrow x)=\infty$.
\end{center}
But the susceptibility is known to be finite at $p_{c}$ for $q>4$ on $\mathbb{Z}^{2}$ (see \cite{DST,DCGHMT}) and is conjectured to be finite for $q>2$ on $\mathbb{Z}^{d}$ with $d\geq 3$ (for $q$ large enough, this result is proved in \cite{laanait1991}). This in turn implies that $p_{c}<\tilde{p}_{c}$ and therefore Corollary 1 is not longer true in these cases.

\medskip

\noindent
2. The argument presented here can be extended to any finite-range coupling constants $(J_{x,y})_{x,y\in\mathbb{Z}^{d}}$, see \cite{HT}.

\medskip

\noindent
3. For infinite-range coupling constants decaying sub-exponentially fast, the second item of Theorem 1 doesn't hold. However, as in \cite{HT}, one can still prove that Lemma 2 holds and that the susceptibility with free boundary conditions is finite for every $p<\tilde{p}_{c}$, which implies $\tilde{p}_{c}=p_{c}$. One can then use the same reasoning as in \cite{Aoun-2020} to deduce that $\mu_{\mathbb{Z}^{d},p,2}(0\leftrightarrow x)\leq cJ_{0,x}$ for every $p<p_{c}$ and for some positive constant $c$ depending on $p$. Note that the use of the exponential decay of the volume of the connected component of $0$ in \cite{Aoun-2020} can be replaced by the existence of $S$ such that $\varphi_ p(S)<1$.
\printbibliography

@article{DST,
   title={Continuity of the Phase Transition for Planar Random-Cluster and Potts Models with {1 $\leq q \leq $4} },
   volume={349},
   number={1},
   journal={Communications in Mathematical Physics},
   publisher={Springer Science and Business Media LLC},
   author={Duminil-Copin, Hugo and Sidoravicius, Vladas and Tassion, Vincent},
   year={2016},
   month={Oct},
   pages={47–107}
}

@book{GRIM2,
author = {Grimmett,G.},
title = {The random-cluster model},
year = {2006},
publisher = {Springer-Verlag},
}

@article{L,
author = {Lieb,E.},
title = {A refinement of Simon's correlation inequality},
journaltitle = {Communications in Mathematical Physics},
year = {1980},
volume = {77},
number = {2},
pages = {127-135},
}

@article{HT,
author = {Duminil-Copin,H. and Tassion,V.},
title = {A new proof of the sharpness of the phase transition for Bernouilli percolation and the Ising model},
journaltitle = {Communications in Mathematical Physics},
year = {2016},
volume = {343},
pages={725–745},
}

@online{DCGHMT,
author = {Duminil-Copin,H. and Gagnebin,M. and Harel,M. and Manolescu,I. and Tassion,V.},
title = {Discontinuity of the phase transition for the planar random-cluster and {P}otts models with {$q >$} 4},
year = {2016},
eprinttype={arXiv},
eprint={1611.09877}, 
}

@article{Aoun-2020,
	title={Sharp asymptotics of correlation functions in the subcritical long-range random-cluster and {P}otts models},
	author={Aoun, Y.},
	year={2020},
	note={arXiv:2007.00116},
}

@article{laanait1991,
author ={Laanait, Lahoussine and Messager, Alain and Miracle-Solé, Salvador and Ruiz, Jean and Shlosman, Senya},
fjournal ={Communications in Mathematical Physics},
journal = {Comm. Math. Phys.},
number = {1},
pages = {81--91},
title ={Interfaces in the {P}otts model. {I}. {P}irogov-{S}inai theory of the {F}ortuin-{K}asteleyn representation},
volume = {140},
year = {1991}
}
\end{document}